\documentclass[11pt]{article}

\usepackage{amsmath,amsfonts,amsthm,amssymb,graphicx,tikz,tikz-cd,url,hyperref,cleveref,mathtools}
\usepackage[all,2cell,ps]{xy}

\theoremstyle{plain}
\newtheorem{thm}{Theorem}[section]

\newtheorem{lem}[thm]{Lemma}
\newtheorem{prop}[thm]{Proposition}

\newtheorem{cor}[thm]{Corollary}

\newtheorem{qtn}[thm]{Question}

\theoremstyle{definition}

\newtheorem{rem}[thm]{Remark}

\theoremstyle{remark}

\newcommand{\bbB}{\mathbb{B}}
\newcommand{\bbC}{\mathbb{C}}

\newcommand{\bbP}{\mathbb{P}}
\newcommand{\bbQ}{\mathbb{Q}}
\newcommand{\bbR}{\mathbb{R}}

\newcommand{\bbV}{\mathbb{V}}

\newcommand{\bbZ}{\mathbb{Z}}

\newcommand{\bfG}{\mathbf{G}}

\newcommand{\bfS}{\mathbf{S}}

\newcommand{\calD}{\mathcal{D}}

\newcommand{\calO}{\mathcal{O}}

\newcommand{\calQ}{\mathcal{Q}}

\newcommand{\calS}{\mathcal{S}}
\newcommand{\calT}{\mathcal{T}}

\newcommand{\al}{\alpha}
\newcommand{\gam}{\gamma}
\newcommand{\Gam}{\Gamma}

\newcommand{\de}{\delta}
\newcommand{\Del}{\Delta}

\newcommand{\lam}{\lambda}
\newcommand{\Lam}{\Lambda}
\newcommand{\sig}{\sigma}
\newcommand{\Sig}{\Sigma}

\DeclareMathOperator{\U}{U}

\DeclareMathOperator{\SL}{SL}
\DeclareMathOperator{\PSL}{PSL}

\DeclareMathOperator{\PU}{PU}
\DeclareMathOperator{\SU}{SU}

\DeclareMathOperator{\id}{id}

\DeclareMathOperator{\Ad}{Ad}
\DeclareMathOperator{\Tr}{Tr}

\DeclareMathOperator{\Isom}{Isom}
\DeclareMathOperator{\Aut}{Aut}

\newcommand{\bs}{\backslash}

\newcommand{\lra}{\longrightarrow}

\newcommand{\conj}{\overline}
\newcommand{\wh}{\widehat}
\newcommand{\wt}{\widetilde}

\newenvironment{pf}{\begin{proof}}{\end{proof}}

\allowdisplaybreaks

\makeatletter
\let\@@pmod\pmod
\DeclareRobustCommand{\pmod}{\@ifstar\@pmods\@@pmod}
\def\@pmods#1{\mkern4mu({\operator@font mod}\mkern 6mu#1)}
\makeatother

\usetikzlibrary{arrows}

\title{Cocompact Fuchsian groups with a modular embedding}
\author{Matthew Stover\\ \small{Temple University}\\ \small{\textsf{mstover@temple.edu}}}
\date{\today}

\begin{document}

\maketitle

\begin{abstract}
A Fuchsian group $\Gam$ has a modular embedding if its adjoint trace field is a totally real number field and every unbounded Galois conjugate $\Gam^\sig$ comes equipped with a holomorphic (or conjugate holomorphic) map ${\phi^\sig : \bbB^1 \to \bbB^1}$ intertwining the actions of $\Gam$ and $\Gam^\sig$ on the Poincar\'e disk $\bbB^1$. This paper provides the first cocompact nonarithmetic Fuchsian groups with a modular embedding that are not commensurable with a triangle group. The main result, proved using period domains, is that any immersed totally geodesic complex curve on a complex hyperbolic $2$-orbifold has a modular embedding. Another consequence is arithmeticity of totally geodesic curves on finite-volume complex hyperbolic surfaces that are commensurable with quotients of $\bbB^1$ by the group generated by reflections in quadrilaterals satisfying certain angle conditions.
\end{abstract}

\section{Introduction}\label{sec:Intro}

This paper considers Fuchsian groups as subgroups of $\PU(1,1)$ acting on the Poincar\'e disk $\bbB^1$. Given such a group $\Gam$, its \emph{adjoint trace field} $\bbQ\!\left(\Tr \Ad(\Gam)\right)$ is the extension of $\bbQ$ generated by the traces of $\Gam$ acting on the Lie algebra of $\PU(1,1)$. A Fuchsian group $\Gam$ then has a \emph{modular embedding} if $\bbQ\!\left(\Tr \Ad(\Gam)\right)$ is a totally real number field and each unbounded Galois conjugate $\Gam^\sig$ induced by a real embedding $\sig$ of $\bbQ\!\left(\Tr \Ad(\Gam)\right)$ comes equipped with a holomorphic (or conjugate holomorphic) map $\phi^\sig : \bbB^1 \to \bbB^1$ intertwining the actions of $\Gam$ and $\Gam^\sig$. In other words,
	\[ \phi^\sig(\gam z) = \gam^\sig \phi^\sig(z) \]
for all $\gam \in \Gam$ and $z \in \bbB^1$. This property is usually stated for subgroups of $\PSL_2(\bbR)$ acting on the upper half-plane (e.g., \cite{CW1, Moller}), but conjugation by the Cayley transformation equates the two definitions.

Arithmetic Fuchsian groups have no unbounded Galois conjugates, so they vacuously have a modular embedding. Cohen and Wolfart proved that triangle groups have modular embeddings \cite{CW1}; see \cite[Thm.\ 6.3]{McMullenHilbert2} on why that paper should allow conjugate holomorphic maps in the definition. The interest in groups with a modular embedding stems from intriguing properties they share with arithmetic Fuchsian groups. For example, their covering radius is a ratio of periods for abelian varieties; these ratios are transcendental, addressing a question of Lang (see \cite{CW1}). They also admit \emph{twisted} modular forms, which were introduced by M\"oller and Zagier \cite{MollerZagier}.

Around a decade after Cohen and Wolfart's work, Calta \cite{Calta} and McMullen \cite{McMullenHilbert} classified Teichm\"uller curves in genus $2$, and McMullen showed that these produce Fuchsian groups with a modular embedding that are not commensurable with triangle groups. See M\"oller \cite{Moller} for examples with higher-degree trace field. It is well-known that Teichm\"uller curves are noncompact, and thus all previously-known nonarithmetic cocompact Fuchsian groups with a modular embedding were subgroups of triangle groups. The first contribution of this paper, proved in \Cref{sec:NotTri}, is the existence of compact examples not commensurable with triangle groups.

\begin{thm}\label{thm:NotTriangle}
There are cocompact nonarithmetic Fuchsian groups with a modular embedding that are not commensurable with any triangle group.
\end{thm}

This gives a negative answer to \cite[Prob.\ 1]{SSW} in the cocompact case. The examples are immersed totally geodesic curves on nonarithmetic complex hyperbolic $2$-orbifolds, i.e., finite-volume quotients of the unit ball $\bbB^2$ in $\bbC^2$ equipped with its Bergman metric of constant holomorphic sectional curvature $-1$ and holomorphic isometry group $\PU(2,1)$. Their modular embedding comes from the following very general result.

\begin{thm}\label{thm:TGModular}
If $\Gam \bs \bbB^1$ is a finite-volume orbifold that immerses as a totally geodesic curve on a finite-volume complex hyperbolic $2$-orbifold, then $\Gam$ has a modular embedding.
\end{thm}

Perhaps surprisingly, the proof does not show that the complex hyperbolic $2$-orbifold has the higher-dimensional analogue of a modular embedding. Instead, it uses the theory of period maps and period domains (more precisely, Mumford--Tate domains) to prove that the Galois conjugate period maps for the lattice in $\PU(2,1)$, sometimes called a \emph{generalized} modular embedding, restrict to give a modular embedding for the Fuchsian subgroup. The proof of \Cref{thm:TGModular}, which relies on the preliminaries in \Cref{sec:Prelims}, is in \Cref{sec:ModularPf}. A higher-dimensional variant of \Cref{thm:TGModular} may hold for totally geodesic curves on complex hyperbolic $n$-orbifolds, but is not explored here in light of the current scarcity of relevant examples.

It is worth noting that the nonarithmetic orbifolds that were constructed by Deligne and Mostow \cite{DM, MostowSigma} may well be inadequate for proving \Cref{thm:NotTriangle}. While Cohen and Wolfart proved that Deligne--Mostow lattices also have the analogue of a modular embedding \cite{CW2}, all known immersed totally geodesic complex curves on nonarithmetic Deligne--Mostow orbifolds are commensurable with triangle orbifolds. In fact, since nonarithmetic complex hyperbolic $2$-orbifolds have only finitely many immersed totally geodesic curves \cite{BFMS2, BU}, it may be the case that all totally geodesic curves on Deligne--Mostow orbifolds are related to triangle orbifolds. The examples in this paper arise from the lattices constructed by Deraux, Parker, and Paupert \cite{DPP2}.

\medskip

\Cref{thm:TGModular} also implies arithmeticity of reflection groups associated with certain properly immersed totally geodesic suborbifolds of finite-volume complex hyperbolic $2$-orbifolds. Let $\calQ[n_1, \dots, n_4]$ be a quadrilateral in $\bbB^1$ with (cyclically ordered) interior angles $\pi / n_j$ and $\wh{\Del}[n_1, \dots, n_4] < \Isom(\bbB^1)$ be the group generated by reflections in the sides of $\calQ[n_1, \dots, n_4]$. Ricker \cite{Ricker}, following up on work of Schmutz Schaller and Wolfart \cite{SSW}, proved for $[n_1, \dots, n_4]$ any of $[2,2,2,t]$, $[2,2,t,t]$, or $[2,t,2,t]$ with $t \in \{3,4,6\}$ that the orientation-preserving index two subgroup $\Del[n_1, \dots, n_4]$ of $\wh{\Del}[n_1, \dots, n_4]$ has a modular embedding if and only if $\wh{\Del}[n_1, \dots, n_4]$ is an arithmetic reflection group. Having a modular embedding is a commensurability invariant, so \Cref{thm:TGModular} and \cite{Ricker} combine to imply the following.

\begin{cor}\label{thm:NoEmbed}
Suppose $[n_1, \dots, n_4]$ is $[2,2,2,t]$, $[2,2,t,t]$, or $[2,t,2,t]$ with $t \in \{3,4,6\}$ and $Z$ is a quotient of $\bbB^1$ commensurable with $\wh{\Del}[n_1, \dots, n_4] \bs \bbB^1$, where $\wh{\Del}[n_1, \dots, n_4]$ is the group generated by reflections in the sides of a quadrilateral of angles $\pi / n_j$. If $Z$ can be holomorphically immersed as a totally geodesic curve on any complex hyperbolic $2$-manifold, then $\wh{\Del}[n_1, \dots, n_4]$ is an arithmetic reflection group.
\end{cor}

This is not the first obstruction to an orbifold being a holomorphically immersed totally geodesic suborbifold of a complex hyperbolic orbifold, but it is of a very different kind than previous obstructions. First, local rigidity of lattices $\Gam < \PU(2,1)$ implies that $\bbQ\!\left(\Tr \Ad(\Gam)\right)$ is a number field. Much less trivially, $\Tr \Ad(\Gam)$ consists only of algebraic integers and $\bbQ\!\left(\Tr \Ad(\Gam)\right)$ is a totally real field \cite[Thm.\ 1.8]{BFMS2}. Thus any Fuchsian group associated with a properly immersed totally geodesic suborbifold of a finite-volume complex hyperbolic $2$-orbifold \emph{semi-arithmetic}, meaning that $\bbQ\!\left(\Tr \Ad(\Gam)\right)$ is a totally real field and $\Tr \Ad(\gam)$ is an algebraic integer for every $\gam \in \Gam$. Thus \Cref{thm:NoEmbed} would follow from known results if any semi-arithmetic group of the appropriate signature were necessarily arithmetic. Schmutz Schaller and Wolfart provide infinitely many distinct signature $(0;2,2,2,t)$ groups that are semi-arithmetic for each $t$ \cite{SSW}, so \Cref{thm:NoEmbed} is indeed new and rules out groups that were unobstructed by previous methods.

\begin{rem}\label{rem:Lots}
A given nonarithmetic group $\Del$ of signature $(0;2,2,2,t)$ can only be commensurable with finitely many groups of signature $(0;2,2,2,s)$, independent of $s$. Indeed, the Margulis commensurability criterion implies that $\Del$ is contained in a unique maximal lattice $\wh{\Del}$ containing all lattices commensurable with $\Del$, and the claim follows easily from the fact that orbifolds of signature $(0;2,2,2,s)$ have uniformly bounded area and $\wh{\Del}$ has finitely many subgroups of bounded finite index. Thus \Cref{thm:NoEmbed} produces infinitely many commensurability classes of examples of semi-arithmetic compact quotients of $\bbB^1$ admitting no proper totally geodesic immersion into a complex hyperbolic $2$-orbifold of finite volume.
\end{rem}

\subsubsection*{Acknowledgements}
Thanks are due to Curtis McMullen for his interest in the paper and for detailed comments on an early draft, Greg Baldi for several conversations related to this paper and the relevant literature, Joshua Lam for noticing a misstatement in an earlier version, and the referee for suggestions. The author acknowledges support from the Institut Henri Poincar\'e (UAR 839 CNRS-Sorbonne Universit\'e) and LabEx CARMIN (ANR-10-LABx-59-01). This material is based upon work supported by Grants DMS-2203555 and DMS-2506896 from the National Science Foundation and the Simons Foundation [SFI-MPS-TSM-00014184, MS].

\section{Preliminaries}\label{sec:Prelims}

\subsubsection*{Basic geometry}

See \cite{Goldman} for more details about basic geometry of $\bbB^2$. It will be advantageous in this paper to consider the model of complex hyperbolic space associated with an arbitrary hermitian form. To that end, if $h$ is any hermitian form of signature $(2,1)$ on $\bbC^3$, the space
	\[ X \coloneqq \left\{\ell \in \bbP^2\ :\ h|_\ell < 0\right\} \]
of $h$-negative lines in projective space $\bbP^2$ can be identified with $\bbB^2$ so that the projective isometry group $G = \PU(h)$, which is isomorphic to $\PU(2,1)$, acts by holomorphic isometries. Let
	\begin{align*}
		K \coloneqq&\ \mathrm{P}\!\left(\U(2) \times \U(1)\right) \\
		\cong&\, \U(2)
	\end{align*}
denote the stabilizer in $G$ of a point in $X$, which determines an identification $X = G / K$ for the left actions of $G$.

Taking $h$-orthogonal complements gives a one-to-one correspondence between $h$-positive lines and totally geodesic embeddings of $\bbB^1$ in $\bbB^2$. Then $\PU(2,1)$ acts transitively on these embeddings with
	\begin{align*}
		W \coloneqq&\ \mathrm{P}\!\left(\U(1, 1) \times \U(1)\right) \\
		\cong&\, \U(1,1)
	\end{align*}
naturally isomorphic to the stabilizer of a fixed Poincar\'e disk $Y \subset X$. Note that the restriction of an isometry to $Y$ is induced by projection of $W$ onto $\PU(1,1)$. The maximal compact subgroup of $W$ is isomorphic to
	\begin{align*}
		L \coloneqq&\ \mathrm{P}\!\left(\U(1) \times \U(1) \times \U(1)\right) \\
		\cong&\, \U(1) \times \U(1)
	\end{align*}
with $W / L$ identified with $\bbB^1$ and one $\U(1)$ factor acting on its normal bundle.

A finite-volume complex hyperbolic $2$-orbifold is then a quotient $\Lam \bs X$ by a lattice $\Lam < G$. Properly holomorphically immersed totally geodesic subspaces
	\[ \Gam \bs \bbB^1 \looparrowright \Lam \bs X \]
will be called \emph{totally geodesic curves} for short. These are determined by $g \in G$ so that
	\[ \Gam_g \coloneqq g W g^{-1} \cap \Lam \]
is a lattice in $g W g^{-1}$, where $\Gam$ is the restriction of $\Gam_g$ to $g(Y) \subset X$. There is then a central exact sequence
	\[ 1 \lra \Del \lra \Gam_g \lra \Gam \lra 1 \]
with $\Del$ the subgroup of $\Gam_g$ acting on the normal bundle to $g(Y)$.

Totally geodesic curves also have a dynamical interpretation. Indeed, they correspond precisely with the closed orbits for the action of $W$ on the homogeneous space $\Lam \bs G$. Specifically, if $[g]$ denotes the point $\Lam g \in \Lam \bs G$, then $[g] W$ is a closed orbit if and only if $\Gam_g$ is a lattice in $g W g^{-1}$. See \cite[\S 8.1]{BFMS2} for more on dynamical characterizations of geodesic subspaces.

\begin{rem}\label{rem:Kobayashi}
A great deal of work has gone into understanding totally geodesic subspaces of locally symmetric manifolds; see for instance \cite{BFMS, BFMS2, BU} and references therein. However, the totally geodesic subvarieties for the \emph{Kobayashi metric} are also of significant interest, especially in the context of quotients of $(\bbB^1)^d$ for $d \ge 2$ and Teichm\"uller space $\calT_g$ with the Teichm\"uller metric; e.g., see \cite{McMullenHilbert, McMullenHilbert2, Moller}. Unlike those other settings, the locally symmetric and Kobayashi metrics coincide for $\bbB^n$.
\end{rem}

\subsubsection*{Arithmetic invariants and period maps}

The next goal is to describe algebraic groups associated with lattices in $\PU(2,1)$ and their associated period maps. See \cite{HodgeBook, GGK} for basic facts about period maps and related objects. Retaining the notation established earlier in this section, let $\Lam < G$ be a lattice and $k = \bbQ\!\left(\Tr \Ad(\Lam)\right)$ denote its adjoint trace field. Then $k$ is totally real and $\Tr(\Ad(\lam))$ is in the ring of integers $\calO_k$ of $k$ for all $\lam \in \Lam$ (see \cite{BFMS2, BU}), but there is in fact more structure.

For all problems considered in this paper, there is no loss of generality in passing to a subgroup of finite index, so $\Lam$ is henceforth assumed to be torsion-free. The adjoint action of $\Lam$ defines a $k$-algebraic group $\bfG$ containing $\Lam$ as a Zariski dense subgroup, and let $\wh{\bfG}$ be the restriction of scalars of $\bfG$ from $k$ to $\bbQ$. After possibly conjugating and passing to a further subgroup of finite index, $\Lam$ is then a subgroup of $\wh{\bfG}(\bbZ)$ where the integral structure is defined through the polarized integral variation of Hodge structure $\wh{\bbV}$ on $\Lam \bs X$ that arises from lifting $\Lam$ to $\SU(2,1)$ and then taking the natural action on $\bbC^3$; see \cite[Thm.\ 1.3.1]{BU}.

Every simple factor of $\wh{\bfG}(\bbR)$ is of Hodge type (e.g., see \cite[Thm.\ 1.5]{BFMS2}), which implies in this case that
	\[ \wh{\bfG}(\bbR) = \prod_{\sig : k \to \bbR} \bfG(k \otimes_\sig \bbR) \cong \PU(2,1)^r \times \PU(3)^s \]
with $r + s = [k : \bbQ]$ and the case $r = 1$ being precisely the arithmetic case. Furthermore, $\wh{\bfG}$ is the derived subgroup of the \emph{generic Mumford--Tate} group of $\wh{\bbV}$ \cite[Cor.\ 5.2.2]{BU}. Associated with $\wh{V}$ is a \emph{period domain} $D = \wh{\bfG}(\bbR) / V$ for $V$ a certain compact subgroup, and $D$ inherits a product structure
	\[ D = X \times \prod_{\sig\, \neq\, \id} D_\sig \]
where $\id$ is associated with the lattice embedding of $\Lam$ and $D_\sig = G_\sig / V_\sig$ for $G_\sig = \bfG(k \otimes_\sig \bbR)$ and $V_\sig$ a compact subgroup of $G_\sig$. The \emph{period map} associated with $\wh{\bbV}$ is then a holomorphic mapping
	\[ \psi : \Lam \bs X \longrightarrow \wh{\bfG}(\bbZ) \bs D \]
that is equivariant for the inclusion of $\Lam$ in $\wh{\bfG}(\bbZ)$. If
	\[ \wt{\psi} : X \longrightarrow X \times \prod_{\sig\, \neq\, \id} D_\sig \]
is the lift to a map on universal coverings, then the first factor is simply the identity map $X \to X$.

As manifolds, the possible target domains $D_\sig$ are in one-to-one correspondence with parabolic subgroups $Q < G_\sig(\bbC) \cong \PSL_3(\bbC)$ so that $Q \cap G_\sig$ is compact and $G_\sig / (Q \cap G_\sig)$ is an open $G_\sig$-orbit in $\PSL_3(\bbC) / Q$. If $Q$ is the stabilizer of a line (or, dually, a plane), then $Q \cap G_\sig$ is conjugate to $K$ and $D_\sig$ is simply $X$ with its unique complex structure (up to complex conjugation). If $Q$ is the minimal parabolic subgroup of $\SL_3(\bbC)$, namely the stabilizer of a full flag, then $V_\sig \cong L$ is the maximal compact subgroup of $W$. This discussion is recorded in the following lemma, where the reader should consult either \cite[p.\ 169-170]{GGK} or \cite[\S 16.4]{HodgeBook} for full details.

\begin{lem}\label{lem:Domains}
The nontrivial period domains for $G = \PU(2,1)$ are
	\begin{align*}
		\bbB^2 &= \PU(2,1) / K \\
		\calD_2 &= \PU(2,1) / L
	\end{align*}
where $K$ is the maximal compact subgroup of $G$ and $L \cong \U(1) \times \U(1)$ is the maximal compact subgroup of the stabilizer $W$ of a totally geodesic $\bbB^1$ in $\bbB^2$.
\end{lem}

Note that there is a projection map $\calD_2 \to \bbB^2$ associated with the natural inclusion of compact subgroups. Specifically, $\calD_2$ is $3$-dimensional and the projection to $\bbB^2$ has fiber $\U(2) / (\U(1) \times \U(1)) \simeq \bbP^1$. In other words, $\calD_2$ is diffeomorphic to $\bbB^2 \times \bbP^1$. While the complex structure on $\bbB^2$ is unique (up to complex conjugation), this is not the case for $\calD_2$, which has several (see \cite[\S IV.G]{GGK} or \cite[p.\ 450-451]{HodgeBook}), and the fibration over $\bbB^2$ may fail to be holomorphic. Specifying a period map $\psi$ will always imply that a complex structure has been fixed on $\calD_2$, but which one will not be relevant in this paper. What will be relevant is the following.

\begin{cor}\label{cor:Subdomains}
Every proper, noncompact Mumford--Tate subdomain of a nontrivial period domain for $\PU(2,1)$ is biholomorphic to $\bbB^1$.
\end{cor}
\begin{pf}
Proper subdomains arise from Hodge type subgroups of $\PU(2,1)$ (e.g., see \cite[\S 15.3]{HodgeBook}), where the only noncompact possibilities up to conjugacy are $W \cong \U(1,1)$ or its derived subgroup, which is isomorphic to $\SU(1,1)$. In both cases, the associated compact subgroup of $\PU(2,1)$ intersects $W$ in its maximal compact subgroup $L$. Indeed, for $\bbB^2$ the compact subgroup of $\PU(2,1)$ is $K$ and $L = W \cap K$, and for $\calD_2$ the compact subgroup is $L$ itself. In particular, the associated subdomain is always biholomorphic to $W / L$, which is $\bbB^1$ (cf.\ \cite[p.\ 177]{GGK}).
\end{pf}

\section{The proof of \Cref{thm:TGModular}}\label{sec:ModularPf}

Retain all notation from \Cref{sec:Prelims}, supposing moreover that $\Lam < G$ is a torsion-free nonarithmetic lattice with $[g] W \subset \Lam \bs G$ the closed $W$-orbit associated with a totally geodesic curve on the complex hyperbolic $2$-manifold $\Lam \bs X$. Let $\Gam$ be the associated finite coarea Fuchsian group, so there is a proper holomorphic totally geodesic immersion
	\[ \Gam \bs \bbB^1 \looparrowright \Lam \bs X \]
and set $\Gam_g = g W g^{-1} \cap \Lam$. Recall that $\Gam_g$ is a central extension of $\Gam$.

The action of $\Gam_g$ on the Lie algebra of $\PU(2,1)$ induces a natural field embedding
	\[ \ell \coloneqq \bbQ\!\left(\Tr \Ad(\Gam)\right) \subseteq k \coloneqq \bbQ\!\left(\Tr \Ad(\Lam)\right) \]
of adjoint trace fields \cite[Rem.\ 5.1.3]{BU}, emphasizing that $\bbQ\!\left(\Tr \Ad(\Gam)\right)$ is the adjoint trace field of $\Gam$ as a lattice in $\PU(1,1)$. Given a real embedding $\sig$ of $\ell$, fix any extension to a real embedding of $k$, which is still denoted by $\sig$. If $\Gam^\sig$ is unbounded, then so is $\Gam_g^\sig < G^\sig$, hence $G^\sig$ is not compact. In other words, unbounded embeddings of $\Gam$ induce unbounded embeddings of $\Lam$.

Associated with $\sig$ is a factor $\wt{\psi}_\sig : X \to D_\sig$ of the lift $\wt{\psi}$ of the period map $\psi$ to $X$. Recall that $g(Y) \subset X$ denotes the totally geodesic embedding of $\bbB^1$ stabilized by $\Gam_g$. A key observation for this paper is the following, which can also be deduced from results in \cite[\S 5]{BU} along with \Cref{cor:Subdomains}, but the proof is not difficult so it is included for completeness.

\begin{prop}\label{prop:KeyProp}
With the notation established in this section and \Cref{sec:Prelims}, the restriction of the period map $\wt{\psi}_\sig$ to $g(Y)$ defines a holomorphic self-map of $\bbB^1$ intertwining the actions of $\Gam$ and $\Gam^\sig$.
\end{prop}
\begin{pf}
Consider the commutative diagram
	\[
		\begin{tikzcd}
			& D_\sig \arrow[dd, "\pi"] \\
			X \arrow[ur, "\wt{\psi}_\sig"] \arrow[dr, "\wt{\rho}_\sig" below left] & \\
			& \bbB^2
		\end{tikzcd}
	\]
where $\pi$ is the natural projection if $D_\sig$ is $\calD_2$ and is the identity if $D_\sig$ is $\bbB^2$. Here $\wt{\rho}_\sig \coloneqq \pi \circ \wt{\psi}_\sig$ intertwines the actions of $\Lam$ and $\Lam^\sig$ on $\bbB^2$, since the period map intertwines the actions and $\pi$ is $G_\sig$-equivariant. Then \Cref{cor:Subdomains} implies that it suffices to prove that $\wt{\psi}_\sig(h(Y))$ is contained in a proper, noncompact Mumford--Tate subdomain of $D_\sig$, since all such subdomains are biholomorphic to $\bbB^1$.

The group $\Gam_g < \Lam$ stabilizing $g(Y)$ fixes the $h$-positive line in $(k \otimes \bbC)^3$ associated with $g(Y) \subset X$. However, extending $\sig$ to an element of $\Aut(\bbC)$, $\Gam_g^\sig < \Lam^\sig$ then also fixes the $\sig$-conjugate line in $(k \otimes_\sig \bbC)^3$, which cannot be a null line. Were the line $h$-positive, then $\Gam_g^\sig$ would fix a point in $\bbB^2$ under the action of $G_\sig \cong \PU(2,1)$, contradicting the assumption that $\Gam_g^\sig$ is unbounded. Therefore $\Gam_g^\sig$ fixes an $h$-positive line, and hence it is contained in a conjugate $g_\sig W g_\sig^{-1}$ of $W$ in $G_\sig$, i.e., $\Gam_g^\sig$ stabilizes a totally geodesic Poincar\'e disk in $\bbB^2$. It follows that $\Gam_g^\sig$ acts on the Mumford--Tate subdomain of $D_\sig$ associated with the orbit $g_\sig W \subset G_\sig$.

It remains to check that $\wt{\psi}_\sig(g(Y))$ is contained in the subdomain associated with $g_\sig W$. This is simply functoriality of period mappings. The restriction of $\wt{\psi}_\sig$ to $g(Y)$ is a direct factor of the period mapping for the restriction of $\wh{\bbV}$ to $\Gam_g$, which by \cite[Thm.\ 5.1.4]{BU} has generic Mumford--Tate group a $\bbQ$-algebraic subgroup of the generic Mumford--Tate group of $\Lam$ that contains the Zariski closure of $\Gam_g^\sig$. The Mumford--Tate domain for $\Gam_g^\sig$ is then a homogeneous space for the Zariski closure of $\Gam_g^\sig$ in $G_\sig$, which is contained in $g_\sig W g_\sig^{-1}$. In particular, the Mumford--Tate domain for $\Gam_g^\sig$ is a proper subdomain of $D_\sig$ containing $g_\sig W / L$, hence it is biholomorphic to the Poincar\'e disk by \Cref{cor:Subdomains}. This completes the proof of the proposition.
\end{pf}

\begin{pf}[Proof of \Cref{thm:TGModular}]
Suppose that $\Gam \bs \bbB^1$ is a finite-volume hyperbolic $2$-orbifold that immerses as a totally geodesic subvariety of a finite-volume complex hyperbolic $2$-orbifold $\Lam \bs \bbB^2$. If $\Lam$ is arithmetic, then it is well-known that $\Gam$ is also arithmetic, so $\Gam$ has a modular embedding. As described earlier in this section, if $\Lam$ is nonarithmetic then $\ell \coloneqq \bbQ\!\left(\Tr \Ad(\Gam)\right)$ is a subfield of $k \coloneqq \bbQ\!\left(\Tr \Ad(\Lam)\right)$. Since $k$ is totally real, so is $\ell$. Any real embedding $\sig$ of $\ell$ then induces an embedding of $k$, and if $\Gam^\sig$ is unbounded, then there is an associated lifted period map $\wt{\psi}_\sig : \bbB^2 \to D_\sig$ for some nontrivial Mumford--Tate domain $D_\sig$ for $\PU(2,1)$. \Cref{prop:KeyProp} implies that $\wt{\psi}_\sig$ provides a holomorphic map $\phi^\sig : \bbB^1 \to \bbB^1$ intertwining the actions of $\Gam$ and $\Gam^\sig$ on $\bbB^1$, and thus $\Gam$ has a modular embedding.
\end{pf}

\begin{rem}\label{rem:Anti}
Note that anti-holomorphic maps are also allowed in the definition of a modular embedding. For the period mapping, the complex structure on the target is chosen so the map is holomorphic.
\end{rem}

\section{The proof of \Cref{thm:NotTriangle}}\label{sec:NotTri}

To prove \Cref{thm:NotTriangle}, it suffices to prove that there is a nonarithmetic compact totally geodesic curve on a nonarithmetic complex hyperbolic $2$-orbifold so that the associated Fuchsian group is not commensurable with a triangle group. The desired example (see \Cref{thm:Others} for others) comes from a nonarithmetic orbifold constructed by Deraux, Parker, and Paupert \cite{DPP2}, where certain special Fuchsian subgroups were classified by Deraux \cite{Deraux}, though the following example is erroneously listed there as having signature $(0;10,10,10,10)$.

\begin{prop}[\S 10 \cite{Deraux}]\label{prop:Martin}
The nonarithmetic Deraux--Parker--Paupert lattice denoted in their notation by $\calS(5, \conj{\sig}_4)$ contains a nonarithmetic Fuchsian group of signature $(0; 2, 10, 10, 10)$ with adjoint trace field $\bbQ(\de)$, where $\de$ is a root of $p(x) = x^4 - 7 x^3 + 9 x^2 + 7 x + 1$.
\end{prop}

Details about this lattice and the Fuchsian subgroup are given for the reader's convenience and to explain the correction in the signature. Define:
	\begin{align*}
		u &\coloneqq e^{2 \pi i / 15} & & & \al &\coloneqq 2 - u^3 - \conj{u}^3 \\
		\tau &\coloneqq -\frac{1+\sqrt{-7}}{2} & & & \beta &\coloneqq (\conj{u}^2 - u)\tau \\
		R_1 &\coloneqq \begin{pmatrix} u^2 & \tau & -u \conj{\tau} \\ 0 & \conj{u} & 0 \\ 0 & 0 & \conj{u} \end{pmatrix} & J &\coloneqq \begin{pmatrix} 0 & 0 & 1 \\ 1 & 0 & 0 \\ 0 & 1 & 0 \end{pmatrix} & H &\coloneqq \begin{pmatrix} \al & \beta & \conj{\beta} \\ \conj{\beta} & \al & \beta \\ \beta & \conj{\beta} & \al \end{pmatrix}
	\end{align*}
Then $\calS(5, \conj{\sig}_4)$ is the image of $\Lam \coloneqq \langle R_1, J \rangle < \SU(H)$ in $\PU(H) \cong \PU(2,1)$. This is a nonarithmetic lattice with adjoint trace field $\bbQ(\de)$ with $\de$ as defined in \Cref{prop:Martin}.

\begin{rem}\label{rem:DerauxChange}
Deraux reports the adjoint trace field of $\calS(5, \conj{\sig}_4)$ as $\bbQ(\beta)$ with $\beta^2 = (5+\sqrt{5})/14$. Here $p(x)$ is the monic integral polynomial generating this extension whose polynomial discriminant equals the discriminant $6125$ of the field. One can check that
	\[ \sqrt{\frac{5+\sqrt{5}}{14}} = \frac{1}{7}\!\left(3 \de^3 - 23 \de^2 + 38 \de + 12\right) \]
under the real embedding of $\bbQ(\de)$ given by
	\[ \de \mapsto \frac{1}{4}\!\left(7 - \sqrt{5} + \sqrt{14(5 - \sqrt{5})}\right)\!. \]
\end{rem}

Note that $R_1$ fixes the first standard basis vector $e_1$, which is $H$-positive and thus defines a totally geodesic embedding $\iota : \bbB^1 \to \bbB^2$ in the model of $\bbB^2$ associated with $H$. Define $R_2 \coloneqq J R_1 J^{-1}$ and $R_3 \coloneqq J^{-1} R_1 J$. One can then directly check that the elements
\begin{align*}
a &\coloneqq (R_1 R_2)^2 & c &\coloneqq R_2 R_3 R_2^{-1} (R_1 J)^2 \\
b &\coloneqq (R_1 R_3)^2 & d &\coloneqq (abc)^{-1}
\end{align*}
all centralize $R_1$ and all have image in $\PU(1,1)$ under restriction to $\iota(\bbB^1)$ of order $10$ except $c$, which has order two since $c^2 = R_1^{-1}$. Thus $\wh{\Gam} \coloneqq \langle a,b,c \rangle$ has restriction to $\iota(\bbB^1)$ a Fuchsian group $\Gam$ of signature $(0; 2,10,10,10)$.

The adjoint trace field of $\Gam$ is also $\bbQ(\de)$. Indeed, the commutator ${\gam = [a,b]}$ of $a$ and $b$ has eigenvalues $\{1,\lam,\lam^{-1}\}$, where $\{\lam, \lam^{-1}\}$ are the eigenvalues associated with the action of $\gam$ on the signature $(1,1)$ subspace associated with $\iota(\bbB^1)$. Thus $(\lam + \lam^{-1})^2$ is $\Tr\Ad(\gam)$ for its action on $\PU(1,1)$ by \cite[Exer.\ 3.3.4]{MacR}. Checking that $\Tr\Ad(\gam)$ has degree four minimal polynomial over $\bbQ$ implies that it must indeed generate $\bbQ(\de)$ over $\bbQ$.

\medskip

The group $\Gam$ is now used to prove \Cref{thm:NotTriangle}.

\begin{pf}[Proof of \Cref{thm:NotTriangle}]
It is noted in \cite{Deraux} that the Fuchsian group $\Gam$ from \Cref{prop:Martin} is nonarithmetic, so it suffices by \Cref{thm:TGModular} to prove that $\Gam$ is not commensurable with a triangle group. Since any lattice in $\PU(1,1)$ containing a triangle group is itself a triangle group \cite[\S 6]{Singerman}, it suffices to show that $\Gam$ is not contained in a triangle group. Let $d$ denote the potential index of $\Gam$ in the triangle group. A simple computer check using orbifold Euler characteristics and the following facts:
\begin{itemize}\itemsep-0.25em
\item[$\star$] some cone point must have order divisible by $10$
\item[$\star$] all cone point orders must be bounded above by $10d$
\item[$\star$] any new torsion must have order divisible by $d$
\end{itemize}
implies the only triangle groups that could contain a group of signature $\Sig \coloneqq (0; 2,10,10,10)$ are given in \Cref{tb:Tris2}, which also records the adjoint trace field of the triangle group, computed by \cite[Exer.\ 4.9.1]{MacR}. One can check by hand that each triangle group contains a subgroup of signature $\Sig$. The adjoint trace field is a commensurability invariant \cite[Thm.\ 3.3.4]{MacR} (and \cite[Exer.\ 3.3.4]{MacR}), and but the invariant trace field of $\Gam$ does not appear in \Cref{tb:Tris2}, so $\Gam$ is not commensurable with a triangle group and the proof of the theorem is complete.
\end{pf}

\begin{table}[h]
\centering
\begin{tabular}{|c|c|}
\hline
Group & Adjoint trace field \\
\hline
(2,3,10) & $\bbQ(\sqrt{5})$ \\
\hline
(3,6,10) & $\bbQ(\cos(\pi / 15))$ \\
\hline
(4,10,20) & $\bbQ(\cos(\pi / 10))$ \\
\hline
\end{tabular}
\caption{Triangle groups containing a group of signature (0;2,10,10,10)}\label{tb:Tris2}
\end{table}

\noindent\textbf{The modular embedding of $\Gam$}

\medskip

Identifying $u^3$ with a primitive $5^{th}$ root of unity $\zeta_5$ identifies the field $\bbQ(\de)$ with a subfield of $\bbQ(\zeta_5, \sqrt{-7})$. Then $\zeta_5 \mapsto \zeta_5^2$ and $\sqrt{-7} \mapsto -\sqrt{-7}$ induce nonidentity Galois automorphisms $\sig_5$ and $\sig_7$ of $\bbQ(\de)$, where the other nontrivial automorphism is induced by $\sig_5 \sig_7$. One can check that the restriction of $H$ to the orthogonal complement to $e_1$ still has signature $(1,1)$ under $\sig_5$ and $\sig_7$ but is definite under $\sig_5 \sig_7$.

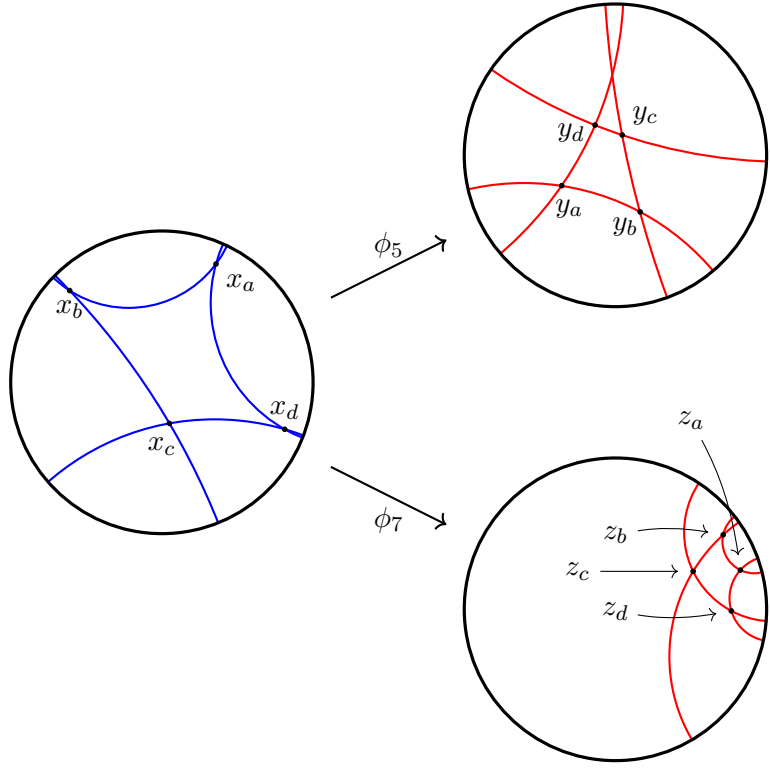
\begin{figure}[h]
\centering
\begin{tikzpicture}[scale=2]
\begin{scope}
\clip (0,0) circle (1cm);
\draw[thick, blue] (-0.22021, 1.21306) circle (0.72109cm);
\end{scope}
\begin{scope}
\clip (0,0) circle (1cm);
\draw[thick, blue] (-4.1434, -2.73332) circle (4.86197cm);
\end{scope}
\begin{scope}
\clip (0,0) circle (1cm);
\draw[thick, blue] (0.34818, -1.91794) circle (1.67324cm);
\end{scope}
\begin{scope}
\clip (0,0) circle (1cm);
\draw[thick, blue] (1.28179, 0.525693) circle (0.958827cm);
\end{scope}
\draw[very thick] (0,0) circle (1cm);
\begin{scope}
\clip (3,1.5) circle (1cm);
\draw[thick, red] (-0.605239+3, -1.82272+1.5) circle (1.6397cm);
\end{scope}
\begin{scope}
\clip (3,1.5) circle (1cm);
\draw[thick, red] (6.92068+3, 1.46829+1.5) circle (7.0037cm);
\end{scope}
\begin{scope}
\clip (3,1.5) circle (1cm);
\draw[thick, red] (1.14055+3, 3.41226+1.5) circle (3.45606cm);
\end{scope}
\begin{scope}
\clip (3,1.5) circle (1cm);
\draw[thick, red] (-2.27036+3, 1.11781+1.5) circle (2.32436cm);
\end{scope}
\draw[very thick] (3,1.5) circle (1cm);
\begin{scope}
\clip (3,-1.5) circle (1cm);
\draw[thick, red] (0.918297+3, 0.448312-1.5) circle (0.210388cm);
\end{scope}
\begin{scope}
\clip (3,-1.5) circle (1cm);
\draw[thick, red] (1.44341+3, -0.311922-1.5) circle (1.08662cm);
\end{scope}
\begin{scope}
\clip (3,-1.5) circle (1cm);
\draw[thick, red] (1.04264+3, 0.509018-1.5) circle (0.588371cm);
\end{scope}
\begin{scope}
\clip (3,-1.5) circle (1cm);
\draw[thick, red] (1.03658+3, 0.0711056-1.5) circle (0.282064cm);
\end{scope}
\draw[very thick] (3,-1.5) circle (1cm);
\draw[->,thick,shorten >=2.5cm,shorten <=2.5cm] (0,0) -- node[above] {$\phi_5$} (3,1.5);
\draw[->,thick,shorten >=2.5cm,shorten <=2.5cm] (0,0) -- node[below] {$\phi_7$} (3,-1.5);
\draw[fill=black] (0.35784,0.78194) node [below right] {$x_a$} circle (0.015cm);
\draw[fill=black] (-0.60985,0.60626) node [below] {$x_b$} circle (0.015cm);
\draw[fill=black] (0.04931,-0.27161) node [xshift=-0.1cm, yshift=-0.3cm] {$x_c$} circle (0.015cm);
\draw[fill=black] (0.81244,-0.3104) node [above] {$x_d$} circle (0.015cm);
\draw[fill=black] (-0.35542+3,-0.20035+1.5) node [xshift=0.1cm, yshift=-0.3cm] {$y_a$} circle (0.015cm);
\draw[fill=black] (0.16351+3,-0.3738+1.5) node [xshift=-0.2cm, yshift=-0.2cm] {$y_b$} circle (0.015cm);
\draw[fill=black] (0.04494+3,0.13445+1.5) node [above right] {$y_c$} circle (0.015cm);
\draw[fill=black] (-0.13458+3,0.20003+1.5) node [yshift=-0.1cm, left] {$y_d$} circle (0.015cm);
\draw[fill=black] (0.82631+3,0.2591-1.5) circle (0.015cm);
\draw[fill=black] (0.71253+3,0.49216-1.5) circle (0.015cm);
\draw[fill=black] (0.51391+3,0.25089-1.5) circle (0.015cm);
\draw[fill=black] (0.76661+3,-0.01059-1.5) circle (0.015cm);
\node at (3.5,-0.25) {$z_a$};
\draw[->,shorten >=0.2cm,shorten <=0.3cm] (3.5,-0.25) to [bend left=10] (0.82631+3,0.2591-1.5);
\node at (3,0.49216-1.5) {$z_b$};
\draw[->,shorten >=0.2cm,shorten <=0.3cm] (3,0.49216-1.5) to [bend left=10] (0.71253+3,0.49216-1.5);
\node at (2.75,0.25089-1.5) {$z_c$};
\draw[->,shorten >=0.2cm,shorten <=0.3cm] (2.75,0.25089-1.5) to (0.51391+3,0.25089-1.5);
\node at (3,-0.01059-1.5) {$z_d$};
\draw[->,shorten >=0.2cm,shorten <=0.3cm] (3,-0.01059-1.5) to [bend right=10] (0.76661+3,-0.01059-1.5);
\end{tikzpicture}
\caption{Quadrilaterals with vertices at the fixed points of $a$, $b$, $c$, and $abc$}\label{fig:Maps}
\end{figure}

This means that $\Gam$ has two nontrivial unbounded embeddings and one bounded embedding, and hence its modular embedding is of the form
	\[ (\id, \phi_5, \phi_7) : \bbB^1 \lra \bbB^1 \times \bbB^1 \times \bbB^1 \]
where $\phi_p : \bbB^1 \to \bbB^1$ intertwines for $\sig_p$, $p \in \{5,7\}$. See \Cref{fig:Maps} for an illustration of $\phi_p$ showing the rectangles in $\bbB^1$ with vertices the fixed points of $a,b,c,d$ under the various embeddings. For analogous images of modular embeddings of triangle groups, see \cite[Fig.\ 6]{McMullenHilbert} and \cite[Fig.\ 5]{McMullenHilbert2}.

\begin{rem}\label{rem:FindMat}
Explicit matrices for the action of $\Gam$ on $\bbB^1$ can easily be found as follows. Set $F = \bbQ(u, \tau)$. Changing coordinates to a basis of $F^3$ consisting of eigenvectors of $c$, which are all defined over $F$, diagonalizes $H$. From there, one can scale the basis to obtain the standard hermitian form on $\bbC^3$ under any embedding of $F$. Under this basis, the matrices for $a,b,c$ become block-diagonal with $2\times 2$ block giving the action on $\bbB^1$. The matrix entries for each are sufficiently complicated as polynomials in $u$ and $\tau$ that they are not included here.
\end{rem}

\begin{rem}\label{rem:NotReflection}
The quadrilateral $\calQ$ with vertices $x_a, x_b, x_c, x_d$ and the neighboring quadrilateral $\calQ^\prime$ with vertices $x_b$, $b^{-1}(x_a)$, $x_c$, and $c(x_d)$ together form a fundamental polygon for the action of $\Gam$. Numerical approximations to the interior angles show that $\calQ$ and $\calQ^\prime$ are not isometric, implying that $\Gam$ cannot be the index two subgroup of a reflection group. In particular, the maps $\phi_p$ cannot be obtained by applying reflections through sides to the Galois conjugate polygons to $\calQ$ shown in \Cref{fig:Maps}. This is in stark contrast with the triangle group setting, where one can apply the Schwarz reflection principle to Galois conjugates of $(p,q,r)$ triangles to build the intertwining maps. The maps $\phi_p$ thus seem much more mysterious than those arising from modular embeddings of triangle orbifolds.
\end{rem}

Lastly, the fact that $\Gam$ has a bounded Galois conjugate means that the associated hyperbolic surface does not naturally immerse in a Hilbert modular $4$-fold. With this and questions from \cite{McMullenHilbert2} in mind:

\begin{qtn}\label{q:Hilbert}
Can totally geodesic curves on nonarithmetic complex hyperbolic manifolds produce new compact totally geodesic curves on Hilbert modular varieties that are not contained in a proper Shimura subvariety?
\end{qtn}

A positive answer would be particularly interesting in dimension $3$. It is known that all compact totally geodesic curves on Hilbert modular surfaces are Shimura curves. McMullen showed that the $(14,21,42)$ triangle group determines a totally geodesic curve on a Hilbert modular $6$-fold not contained in a proper Shimura subvariety \cite[Thm.\ 1.9]{McMullenHilbert2}.


\medskip

\noindent\textbf{Other examples}

\medskip

The example presented in this paper is worked out in detail because it is the simplest among a trio of candidates from the lattices $\calS(p, \conj{\sig}_4)$ in \cite{DPP2}. Using the methods in this section for the lattices $\calS(8, \conj{\sig}_4)$ and $\calS(12, \conj{\sig}_4)$, one can similarly obtain the following result, where the examples are obtained with literally no change to the definitions for the previous example except changing the value $e^{2 \pi i / 15}$ to $e^{2 \pi i / 3 p}$. (Note that the signatures stated in \cite{Deraux} are again erroneous, for the same reason as before that $c^2 = R_1$.)

\begin{thm}\label{thm:Others}
There are nonarithmetic Fuchsian groups $\Gam_1$ and $\Gam_2$ with respective signatures $(0;2,4,4,8)$ and $(0;2,3,3,4)$ that have a modular embedding and are not commensurable with each other, any triangle group, or the example of signature $(0;2,10,10,10)$ from this section.
\end{thm}
\begin{pf}[Sketch of the proof]
The lattice $\calS(8, \conj{\sig}_4)$ contains a group $\Gam_1$ of signature $(0;2,4,4,8)$ with adjoint trace field $\bbQ(\sqrt{2}, \sqrt{7})$. The triangle groups containing a subgroup isomorphic to $\Gam_1$ are $(2,4,32)$, $(2,5,8)$, $(2,8,12)$, $(3,4,8)$, and $(4,4,16)$, none of which have adjoint trace field isomorphic to $\bbQ(\sqrt{2}, \sqrt{7})$, hence $\Gam_1$ is not commensurable with a triangle group. The group $\calS(12, \conj{\sig}_4)$ similarly contains a group $\Gam_2$ of signature $(0;2,3,3,4)$ with adjoint trace field $\bbQ(\sqrt{3}, \sqrt{7})$, where the relevant triangle groups $(2,3,8)$, $(2,3,12)$, $(2,3,20)$, $(3,3,4)$, and $(3,4,8)$ are ruled out similarly. Both $\Gam_1$ and $\Gam_2$ have a modular embedding by \Cref{thm:TGModular}, and Deraux \cite{Deraux} checked that they are nonarithmetic. All the examples in the section have distinct adjoint trace fields, so they are mutually incommensurable. This completes the sketch.
\end{pf}

\begin{rem}\label{rem:Alas}
Each example from \Cref{thm:Others} has a bounded embedding, hence also does not answer \Cref{q:Hilbert} in the positive.
\end{rem}

\begin{rem}\label{rem:Almost}
Another possible example is from the lattice $\calT(5,\bfS_2)$ in \cite{DPP2}, which may be contained in a triangle group. Either way, from \cite{CW1} or \Cref{thm:TGModular}, the group has a modular embedding.
\end{rem}

\bibliography{ModularEmbeds}

\end{document}